\documentclass{amsart}

\usepackage[T1]{fontenc}
\usepackage[utf8]{inputenc}
\usepackage{textcomp}
\usepackage{mathtools}
\usepackage{amstext}
\usepackage{amsthm}
\usepackage{amssymb}
\usepackage{graphicx}
\usepackage{hyperref}
\makeatletter

\providecommand{\tabularnewline}{\\}
\newcommand{\lyxdot}{.}

\numberwithin{equation}{section}
\numberwithin{figure}{section}
\theoremstyle{plain}
\newtheorem{thm}{\protect\theoremname}
\theoremstyle{plain}
\newtheorem{lem}[thm]{\protect\lemmaname}
\theoremstyle{plain}
\newtheorem{cor}[thm]{\protect\corollaryname}
\theoremstyle{remark}
\newtheorem{rem}[thm]{\protect\remarkname}
\theoremstyle{plain}
\newtheorem{algorithm}[thm]{\protect\algorithmname}
\newenvironment{lyxlist}[1]
	{\begin{list}{}
		{\settowidth{\labelwidth}{#1}
		 \setlength{\leftmargin}{\labelwidth}
		 \addtolength{\leftmargin}{\labelsep}
		 }}
	{\end{list}}
\theoremstyle{definition}
\newtheorem{example}[thm]{\protect\examplename}

\newcommand{\CK}{\ensuremath{ \mathbb{C}_{K} } }
\newcommand{\RK}{\ensuremath{\mathbb{R}_{K}} }
\renewcommand{\H}{\ensuremath{\mathbb{H}} }
\newcommand{\T}{\ensuremath{\mathbb{T}} }
\newcommand{\M}{\ensuremath{\mathbb{M}} }

\newcommand{\HK}{\ensuremath{\mathbb{H}_{K}} }
\newcommand{\OK}{\ensuremath{\mathcal{O}_{K}} }
\newcommand{\U}{\ensuremath{\mathcal{U}} }
\newcommand{\PGL}{\ensuremath{\mathrm{PGL}}}

\newcommand{\PSLZ}{\ensuremath{\mathrm{PSL}_{2}(\mathbb{Z}) } }
\newcommand{\PSLK}{\ensuremath{\mathrm{PSL}_{2}(K) } }

\newcommand{\SLOK}{\ensuremath{\mathrm{SL}_{2}(\mathcal{O}_{K}) } }

\newcommand{\GLR}{\ensuremath{\mathrm{GL}_{2}(\mathbb{R}) } }
\newcommand{\PSLOK}{\ensuremath{\mathrm{PSL}_{2}(\mathcal{O}_{K})} }

\newcommand{\ida}{\ensuremath{\mathfrak{a}} }
\newcommand{\idb}{\ensuremath{\mathfrak{b}} }
\newcommand{\idc}{\ensuremath{\mathfrak{c}} }
\renewcommand{\P}[1]{\ensuremath{\mathbb{P}^{1}(#1)} }
\newcommand{\PK}{\ensuremath{\P{K}} }
\newcommand{\Tr}{\ensuremath{\mathrm{Tr}} }
\newcommand{\G}{\ensuremath{\mathcal{C}} }

\let\textv\v
\renewcommand{\vec}[1]{\mathbf{#1}}
\renewcommand\v{\TextOrMath{\textv}{\vec}}
\newcommand{\Hn}{\ensuremath{\mathbb{H}^{n}} }

\newcommand{\N}{\ensuremath{\mathrm{N}} }

\newcommand{\emb}{\ensuremath{\varphi}}

\newcommand{\cusp}[2]{\ensuremath{(#1:#2)}}
\usepackage{tikz,ifthen}
\usepackage{subfig}
\usetikzlibrary{calc,arrows}
\usepackage{mathtools}
\usepackage{etoolbox}
\usetikzlibrary{decorations.markings,patterns,decorations.pathreplacing,intersections}
\usetikzlibrary{backgrounds}
 \usepackage{pgf}
\DeclareUnicodeCharacter{2212}{-}

\makeatother

\usepackage{babel}
\addto\captionsamerican{\renewcommand{\algorithmname}{Algorithm}}
\addto\captionsamerican{\renewcommand{\corollaryname}{Corollary}}
\addto\captionsamerican{\renewcommand{\examplename}{Example}}
\addto\captionsamerican{\renewcommand{\lemmaname}{Lemma}}
\addto\captionsamerican{\renewcommand{\remarkname}{Remark}}
\addto\captionsamerican{\renewcommand{\theoremname}{Theorem}}
\addto\captionsenglish{\renewcommand{\algorithmname}{Algorithm}}
\addto\captionsenglish{\renewcommand{\corollaryname}{Corollary}}
\addto\captionsenglish{\renewcommand{\examplename}{Example}}
\addto\captionsenglish{\renewcommand{\lemmaname}{Lemma}}
\addto\captionsenglish{\renewcommand{\remarkname}{Remark}}
\addto\captionsenglish{\renewcommand{\theoremname}{Theorem}}
\providecommand{\algorithmname}{Algorithm}
\providecommand{\corollaryname}{Corollary}
\providecommand{\examplename}{Example}
\providecommand{\lemmaname}{Lemma}
\providecommand{\remarkname}{Remark}
\providecommand{\theoremname}{Theorem}

\title{A reduction algorithm for Hilbert modular groups}

\author{Fredrik Str\"omberg}
\thanks{This work was partially supported by the Engineering and Physical Sciences Research Council [EP/N007360/1]}
\address{
School of Mathematical Sciences,
The University of Nottingham,
University Park, \\
Nottingham NG7~2RD,
United Kingdom
}
\email{fredrik.stromberg@nottingham.ac.uk}

\begin{document}

\begin{abstract}
The aim of this paper is to present an explicit reduction algorithm for Hilbert modular groups over arbitrary totally real number fields.
An implementation of the algorithm is available to download from \cite{hilbertmodgrp}.  The exposition is self-contained 
 and sufficient details are given for the reader to understand how it works and implement their own version if desired.
\end{abstract}

\maketitle

\section{Introduction }

Given a group $G$,  acting on a topological space $X$, it is often
useful to have a set of representatives of the orbit,  $G\backslash X$,
which are ``reduced'' with respect to some suitable definition.
In number theory the most prominent example is the reduction theory
of the modular group $\Gamma=\PSLZ$. This gives rise to a large number
of interesting applications including the classical theory of Gauss
for binary quadratic forms and continued fractions, as well as more
recent developments in modular and automorphic forms. 

In the case of the modular group the topological space can be viewed
as the complex upper-half plane $\H=\{z=x+iy\in\mathbb{C}\mid y>0\}$
and the action is given by Möbius transformations. The most commonly
used set of representatives of $\Gamma\backslash\H$ is given by the
following set
\[
\mathcal{F}=\left\{ z=x+iy\in\mathbb{C}\mid-1/2\le x\le1/2,\ \left|z\right|\ge1\right\} .
\]
This is an example of a \emph{closed fundamental domain,} meaning
that it tessellates the upper half-plane, $\H=\Gamma\mathcal{F}$,
and different copies overlap only on the boundary, i.e. $V\mathcal{F}^{\circ}\cap W\mathcal{F}^{\circ}=\emptyset$
if $V\neq W\in\Gamma.$ By the covering property it is clear that
for any $z\in\H$ there exists some $A\in\Gamma$ such that $Az\in\mathcal{F}$
and it is easy to see that unless $z$ is a fixed point of $\Gamma$
this element $A$ is unique. This geometric reduction can then be
translated into a reduction theory of, for instance, binary quadratic
forms, by noting that the action of $\Gamma$ on $q(x,y)=ax^{2}+bxy+cy^{2}$
with discriminant $\Delta=b^{2}-4ac<0$ is equivalent to the action
of $\Gamma$ on the point $x_{0}=\frac{1}{2a}(-b+i\sqrt{\left|\Delta\right|})$
in the upper half-plane. 

The goal of the current paper is to give a coherent presentation of
an explicit algorithm which \foreignlanguage{american}{generalizes}
the reduction procedure to the case of Hilbert modular groups, $\Gamma_{K}=\PSLOK,$
for arbitrary totally real number fields $K$. With notation as in
the following sections our main result is the following. 
\begin{thm}
Given a totally real number field $K$ and $\vec{z}=\vec{x}+i\vec{y}\in\HK$
there exists an explicit algorithm (Algorithm \eqref{alg:Reduction-algorithm}),
which finds an element $A\in\Gamma_{K}$ such that $A\vec{z}\in\mathcal{F}_{K},$
where $\mathcal{F}_{K}$ is a certain fundamental domain for $\Gamma_{K}$.
Furthermore, the runtime of this algorithm is polynomial in $\N(\vec{x}),$
$\N(\vec{y}),$ $\N$$\left(\vec{y}^{-1}\right)$ as $\vec{z}$ varies.
\end{thm}

We will follow a construction of fundamental domains for Hilbert modular
groups that originated with Blumenthal \cite{MR1511187} and was further
improved by e.g. Maaß \cite{MR0003405}, Herrmann \cite{MR0062773}
and Tamagawa \cite{MR0121356}. A comprehensive description of this
method is given in the lecture notes by Siegel \cite[Ch.~3.2]{MR659851}
and it is this presentation we have mainly followed. 

The fundamental domains for quadratic fields and in particular $\mathbb{Q}(\sqrt{5})$
and some others of class number one have been studied in more detail
both theoretically and numerically by e.g. G\"otzky \cite{MR1512493},
Cohn \cite{MR195818,MR0174528}, Deutsch \cite{MR2575559,MR1898756},
Jespers, Kiefer and del R\'{\i}o \cite{MR3511291}, and Quinn and
Verjovsky \cite{MR4091535}.

The intention with the presentation of this paper is to make sure
that the exposition is as self-contained as possible and that all
details and notations are clear to the reader, in particular, regarding
which groups we consider. The following three chapters, 2, 3 and 4,
are therefore mainly aiming at reformulating elementary results from
mainly Siegel \cite[Ch.~3.2]{MR659851} and van der Geer \cite{MR930101}
but also other sources, into a common language. We start with a brief
summary of number fields and embeddings, followed by a section on
Hilbert modular groups and the different types of elements. After
this we give a theoretical presentation of the fundamental domain
and the different components involved. This is followed by Section
Section 5, where a detailed analysis of the proof of the existence
of a closest cusp gives rise to Theorem \eqref{thm:sigma_rho_bound_emb},
which is the theoretical foundation behind the algorithm. 

After all necessary theoretical results are presented we then give
the actual reduction Algorithm, separated in two algorithms to be
more comprehensive. After this we provide a selection of detailed
examples with the aim to demonstrate the veracity and effectiveness
of the algorithm, covering examples of class number greater than one
and degree greater than two, both of which previous numerical methods
have not been able to deal with successfully. As a conclusion we mention
some proposed further work and applications. 

It should be noted that all algorithms mentioned in this paper are
implemented using SageMath \cite{sage} and is available as a Python
package at \cite{hilbertmodgrp}. Furthermore, all examples presented
in Section \eqref{sec:Examples} (and more) are available in Jupyter
notebook format as part of this package. 

\section{Number Fields and Embeddings}

Let $K$ be a totally real number field of degree $n$ over $\mathbb{Q}$
with ring of integers $\OK$ and unit group $\U$. Choose an integral
basis $\alpha_{1},\ldots,\alpha_{n}$ of $\OK$ and a set of generators
(fundamental units) $\varepsilon_{1},\ldots,\varepsilon_{n-1}$ of
$\U$. Let $\emb{}_{i}:K\hookrightarrow\mathbb{R},$ $i=1,\ldots,n$
be the embeddings of $K$ into $\mathbb{R}$ and define the norm and
trace on $K/\mathbb{Q}$ by
\[
\N\coloneqq\N_{K/\mathbb{Q}}:\alpha\mapsto\prod\emb_{i}\left(\alpha\right)\quad\text{and}\quad\Tr\coloneqq\Tr_{K/\mathbb{Q}}:\alpha\mapsto\sum\emb_{i}\left(\alpha\right).
\]
When there is no risk of confusion we sometimes write $\alpha_{i}$
for $\emb_{i}(\alpha$). The ideal class number of $K$ is denoted
by $h$ and we let $\idc_{1},\ldots,\idc_{h}$ be the set of ideal
classes, with $\idc_{1}$ the trivial class, and $\ida_{1}=(1),\ida_{2},\ldots,\ida_{h}$
a fixed set of ideal class representatives, chosen by selecting a
fixed ideal of smallest norm in each class.

An element $\alpha\in K$ is said to be totally positive, and we write
$\alpha\gg0$, if $\emb_{i}(\alpha)>0$ for all embeddings $\emb_{i}$.
To further simplify certain formulas we introduce the rings $\CK=\mathbb{C}\otimes_{\mathbb{Q}}K$
and $\RK=\mathbb{R}\otimes_{\mathbb{Q}}K$ and view $\CK$ as an algebra
over both $\mathbb{C}$ and $K$ with the multiplication operations
defined in the natural way. More precisely, for pure tensors $\vec{z},\vec{z'}\in\CK$
with $\vec{z}=z\otimes a$ and $\vec{z}'=z'\otimes a'$ for some $z,z'\in\mathbb{C}$
and $a,a'\in K$ we define 
\begin{align*}
\vec{z}\vec{z}' & =zz'\otimes aa',\quad z'\vec{z}=\vec{z}z'=(z'z)\otimes a,\quad a'\vec{z}=\vec{z}a'=z\otimes(a'a),
\end{align*}
and then extend these operations to the whole of $\CK$ by linearity,
and similarly for elements of $\RK$. The real and imaginary parts
of $\vec{z}=z\otimes a$ are defined by
\[
\Im(\vec{z})=\Im(z)\otimes a\in\RK\quad\text{and}\quad\Re(\vec{z})=\Re(z)\otimes a\in\RK,
\]
again extended linearly, and we will write a general $\vec{z}\in\mathbb{C}_{K}$
as $\vec{z}=\vec{x}+i\vec{y}$ with $\vec{x}=\Re(\vec{z})$ and $\vec{y}=\Im(\vec{z})$.
The embeddings $\emb_{i}$ are extended to embeddings of $\CK$ in
$\mathbb{C}$ and $\RK$ in $\mathbb{R},$ respectively, by setting
\[
\emb_{i}(\vec{z})=\emb_{i}(z\otimes a)=z\emb_{i}(a)
\]
and we use these to define the trace and norm on $\CK$ and $\RK$.
An element $\vec{x}\in\RK$ is said to be totally positive, written
$\vec{x}\gg0$, if $\emb_{i}(\vec{x})>0$ for all embeddings $\emb_{i}$
and similarly we write $\vec{x}\gg\vec{y}$, or equivalently, $\vec{y}\ll\vec{x}$,
if $\vec{x}-\vec{y}\gg0.$ If $\vec{z}\in\mathbb{C}_{K}$ then $\left|\vec{z}\right|\in\mathbb{R}_{K}$
is defined by $\emb_{i}(\left|\vec{z}\right|)=\left|\emb_{i}(\vec{z})\right|$
for all $i$. We define an analog of the standard upper half-plane
by setting 
\[
\HK=\left\{ \vec{z}\in\CK\mid\Im(\vec{z})\gg0\right\} .
\]
Many classical results about Hilbert modular groups and forms are
formulated in terms of $n$ copies of the standard upper half-plane
\[
\H^{n}=\left\{ (z_{1},\ldots,z_{n})\in\mathbb{C}^{n}\mid\Im(z_{i})>0,\,1\le i\le n\right\} 
\]
but it is very easy to translate results between this and $\mathbb{H}_{K}$
using the embedding $\emb$ of $\CK$ into $\mathbb{\mathbb{C}}^{n}$
(as vector spaces) given by
\[
\vec{z}\mapsto\emb\left(\vec{z}\right)=\left(\emb_{1}(\vec{z}),\ldots,\emb_{n}(\vec{z})\right)\in\mathbb{C}^{n}.
\]

\section{Hilbert Modular Groups}

For the purpose of this paper it is most natural to define the Hilbert
modular group for $K$ as the projective group

\begin{align*}
\Gamma_{K} & =\PSLOK\simeq\SLOK/\{\pm I_{2}\},
\end{align*}
where $I_{2}$ is the $2$-by-$2$ identity matrix, and we usually
represent the elements of $\Gamma_{K}$ by the associated matrices.
In connection with cusps it is also natural to consider the following
group associated with an integral ideal $\idb$ of $K$:
\[
\Gamma\left(\OK\oplus\idb\right)=\left\{ \left(\begin{smallmatrix}\alpha & \beta\\
\gamma & \delta
\end{smallmatrix}\right),\,\alpha,\delta\in\OK,\beta\in\idb^{-1},\gamma\in\idb,\,\alpha\delta-\beta\gamma=1\right\} \subseteq\PSLK.
\]
 The group $\PSLK$ acts on $\HK$ by linear fractional transformations:

\begin{equation}
A(\vec{z})=\frac{\alpha\vec{z}+\beta}{\gamma\vec{z}+\delta}\coloneqq(\alpha\vec{z}+\beta)\left(\gamma\vec{z}+\delta\right)^{-1}\quad\text{if }A=\left(\begin{smallmatrix}\alpha & \beta\\
\gamma & \delta
\end{smallmatrix}\right)\in\PSLK,\label{eq:Hn_action}
\end{equation}
and this action is extended as usual to $\PK$ by setting 
\begin{equation}
A\left(\rho:\sigma\right)=(\alpha\rho+\beta\sigma:\gamma\rho+\delta\sigma)\quad\text{if}\quad(\rho:\sigma)\in\PK.\label{eq:PK_action}
\end{equation}
Elements of $\PSLK$ can be classified, for instance, by using the
trace of the associated matrix. For convenience we use the same terminology
as in $\text{GL}_{2}(\mathbb{R})$ and we say that $A$ is:
\begin{itemize}
\item \emph{parabolic} if $\Tr(A)=\pm2$,
\item \emph{elliptic} if $\left|\Tr\left(A\right)\right|\ll2$, and
\item \emph{hyperbolic} if $\left|\Tr(A)\right|\gg2$.
\end{itemize}
It is clear that $A$ is elliptic, parabolic or hyperbolic precisely
if all embeddings $\emb_{i}(A)$ are of the corresponding type in
$\GLR.$ An element that does not belong to any of these types is
simply said to be \emph{mixed.} It is not hard to show that $A$ is
parabolic if and only if it has a unique fixed point in $\PK,$ elliptic
if and only if it has a unique fixed point in $\HK$ and hyperbolic
if and only if it has two fixed points in $\PK$. For more details
see e.g. Freitag \cite[II.§2-§3]{MR1050763}

\subsection{Element and generators of the Hilbert modular group}

If $\alpha\in\OK$ and $\varepsilon\in\U$ we define the following
elements of $\Gamma_{K}$: 
\begin{align*}
T^{\alpha} & \coloneqq\left(\begin{array}{cc}
1 & \alpha\\
0 & 1
\end{array}\right),\quad E(\varepsilon)\coloneqq\left(\begin{array}{cc}
\varepsilon & 0\\
0 & \varepsilon^{-1}
\end{array}\right)\quad\text{and}\ S\coloneqq\left(\begin{array}{cc}
0 & -1\\
1 & 0
\end{array}\right).
\end{align*}
and note that the corresponding actions on $\mathbb{H}_{K}$ are given
by the maps
\[
T^{\alpha}:\vec{z}\mapsto\vec{z}+\alpha,\quad E(\varepsilon):\vec{z}\rightarrow\varepsilon^{2}\vec{z}\quad\text{and}\ S:\vec{z}\rightarrow-\vec{z}^{-1}.
\]
For an integral ideal $\ida\subseteq\OK$ we define the \emph{translation
module} of $\ida$ by 
\[
\T^{\ida}\coloneqq\left\{ T^{\beta}\mid\beta\in\ida\right\} 
\]
and if $\mathcal{H}\leqslant\U$ is generated by $\varepsilon_{1}^{a_{1}},\varepsilon_{2}^{a_{2}},\ldots,\varepsilon_{n-1}^{a_{n-1}}$
then the set of\emph{ multipliers} of $\mathcal{H}$ is
\[
\M_{\mathcal{H}}\coloneqq\left\{ E(\varepsilon)\mid\varepsilon\in\mathcal{H}\right\} \simeq\left\langle E(\varepsilon_{1}^{a_{1}})\right\rangle \times\cdots\times\left\langle E(\varepsilon_{n-1}^{a_{n-1}})\right\rangle .
\]
Let $\beta_{1},\ldots,\beta_{n}$ and $\beta_{1}',\ldots,\beta_{n}'$
be integral bases of $\ida$ and $\ida^{-1}$. It is clear that the
translation modules are finitely generated, more precisely
\[
\T^{\ida}\simeq\left\langle T^{\beta_{1}}\right\rangle \times\cdots\times\left\langle T^{\beta_{n}}\right\rangle \quad\text{and}\quad\T^{\ida^{-1}}\simeq\left\langle \smash{T^{\beta_{1}'}}\vphantom{T^{\beta}}\right\rangle \times\cdots\times\left\langle \smash{T^{\beta_{n}'}}\vphantom{T^{\beta}}\right\rangle .
\]
It follows by a result of Vaser\v{s}te\u{\i}n \cite{MR0435293} (see
also \cite{MR611449} and \cite[p.~82]{MR930101}) that $\Gamma(\OK\oplus\ida)$
is generated by upper and lower-triangular matrices and since these
can all be expressed in terms of the elements $S$ and $T^{\alpha}$
it is generated by the set
\[
\{S,T^{\beta_{1}},T^{\beta_{2}},\ldots,T^{\beta_{n}},T^{\beta_{1}'},\ldots,T^{\beta_{n}'}\}.
\]
As a special case we conclude that $\Gamma_{K}=\Gamma(\mathcal{O}_{K}\oplus\mathcal{O}_{K})$
is generated by $\{S,T^{\alpha_{1}},\ldots,T^{\alpha_{n}}\}$. This
set of generators is very simple and an immediate extension of the
well-known generators $S$ and $T=T^{1}$ for $\PSLZ$. Unfortunately,
in the case of $\Gamma_{K}$, these generators do not have the same
geometric significance and in particular do not correspond to side-pairing
transformations. They are therefore not immediately useful in a reduction
algorithm. It is therefore common to consider a slightly larger set
of generators including elements $E(\varepsilon)$ with $\varepsilon\in\U$
even though these can of course be expressed by the other generators
using, for instance, the algorithms introduced in \cite{MR3290967}. 

\subsection{Cusps of Hilbert Modular Groups\label{subsec:Cusps-of-Hilbert}}

The set of \emph{cusps} of $\Gamma_{K}$, in other words, fixed points
of parabolic elements, can be identified with the projective line
$\PK=K\cup\left\{ \infty\right\} $ where the cusp at infinity, $\infty$,
is as usual a convenient symbol for the class $(1:0)$. Most results
below are well-known and for proofs and further details we refer the
reader to e.g. \cite{MR930101} or \cite{MR3381450}. Note that $\Delta(\OK,\idb^{-1})$
in the notation of \cite{MR3381450} corresponds to the group $\Gamma(\mathcal{O}_{K}\oplus\mathfrak{b})$
in our notation.

Every cusp $\lambda\in\PK$ can be represented by $\cusp{\rho}{\sigma}$
for some non-unique pair $\rho,\sigma\in\OK$ with associated fractional
ideal $\ida_{\rho,\sigma}=\left(\rho,\sigma\right)$. It is easy to
see that different representatives for $\lambda$ give rise to fractional
ideals in the same ideal class, denoted by $\idc_{\lambda}\in\text{Cl}(K)$.
For any $\lambda\in\PK$ we assume that $\rho$ and $\sigma$ are
chosen such that $(\rho,\sigma)=\ida_{j}$ for some ideal class representative
$\ida_{j}$. 

Furthermore, the ideals associated with $\left(\rho:\sigma\right)$
and $A\left(\rho:\sigma\right)$ are identical if $A\in\Gamma_{K}$
since $\det(A)=1$. It can be shown that the map $\lambda\mapsto\idc_{\lambda}$
is a bijection from $\Gamma_{K}\backslash\PK\rightarrow\text{Cl}(K)$
and therefore the number of $\Gamma_{K}$-equivalence classes of cusps
is equal to $h$, the ideal class number of $K$, and we choose $\lambda_{1}=\infty,\ldots,\lambda_{h}$
as representatives for $\Gamma\backslash\PK$ such that $\lambda_{j}$
is associated with $\idc_{j}$ and we write $\lambda_{j}=\cusp{\rho_{j}}{\sigma_{j}}$
with $(\rho_{j},\sigma_{j})=\ida_{j}$. 

It is easy to see that the stabilizer of the cusp $\infty$ in $\Gamma_{K}$
is given by 
\begin{eqnarray*}
\Gamma_{K,\infty} & :=\left\{ T^{\alpha}E(\varepsilon)=\left(\begin{smallmatrix}\varepsilon & \varepsilon^{-1}\alpha\\
0 & \varepsilon^{-1}
\end{smallmatrix}\right):\vec{z}\mapsto\varepsilon^{2}\vec{z}+\alpha,\;\varepsilon\in\U,\;\alpha\in\OK\right\} \simeq & T^{\OK}\rtimes\U^{2}.
\end{eqnarray*}
Corresponding to each cusp representative $\lambda_{j}=(\rho_{j}:\sigma_{j})$
we choose a cusp normalizing map, $A_{j}\in\PGL_{2}(\mathcal{O}_{K})$,
such that $A_{j}\left(\infty\right)=\lambda_{j}$ and 
\[
A_{j}=\left(\begin{array}{cc}
\rho_{j} & \xi_{j}\\
\sigma_{j} & \eta_{j}
\end{array}\right)
\]
with $\xi_{j},\eta_{j}\in\ida_{j}^{-1}$ and $\rho_{j}\eta_{j}-\sigma_{j}\xi_{j}=1$.
In the notation of \cite{MR3381450} $A_{j}$ is an $(\ida_{j},\ida_{j}^{-1})$-matrix.
The map $A_{j}$ is unique up to multiplication by an element in $\Gamma_{K,\infty}$
on the right and we have
\[
A_{j}^{-1}\Gamma_{K}A_{j}=\Gamma(\OK\oplus\ida_{j}^{2}).
\]
As an alternative to studying the set of cusp representatives $\lambda_{1},\ldots,\lambda_{h}$
of $\Gamma_{K}$ it is therefore possible to consider the cusp at
$\infty$ for the collection of groups $\Gamma\left(\OK\oplus\ida_{j}^{2}\right)$
for $j=1,\ldots,h$, with stabilizers 
\begin{eqnarray*}
\Gamma(\OK\oplus\ida_{j}^{2})_{\infty} & = & \left\{ T^{\alpha}E(\varepsilon)=\left(\begin{smallmatrix}\varepsilon & \varepsilon^{-1}\alpha\\
0 & \varepsilon^{-1}
\end{smallmatrix}\right):\vec{z}\mapsto\varepsilon^{2}\vec{z}+\alpha,\;\varepsilon\in\U,\;\alpha\in\ida_{j}^{-2}\right\} \\
 & \simeq & T^{\ida_{j}^{-2}}\rtimes\U^{2}.
\end{eqnarray*}
For an arbitrary cusp $\mu\in\PK$ we choose a map $U_{\mu}\in\Gamma_{K}$
such that $U_{\mu}\left(\mu\right)=\lambda_{I\left(\mu\right)}$ where
$I(\mu)$ is a unique integer in $\{1,2,\ldots,h\}$ and define the
cusp normalizer of $\mu$ as
\[
A_{\mu}=U_{\mu}^{-1}A_{I\left(\mu\right)}.
\]
 It is now easy to show (see e.g. \cite{MR659851}) that the stabilizer
of an arbitrary cusp $\mu\in K$ in $\Gamma_{K}$ can be written as
\[
\Gamma_{K,\mu}=A_{\mu}\Gamma(\OK\oplus\ida_{I(\mu)}^{2})_{\infty}A_{\mu}^{-1}
\]
where $\ida_{I(\mu)}$ is the ideal corresponding to the cusp representative
$\lambda_{I(\mu)}$ which is equivalent to $\mu$. In particular,
all elements that stabilizes $\mu$ in $\Gamma_{K}$ can be written
as $A_{\mu}T^{\alpha}E(\varepsilon)A_{\mu}^{-1}$ for some $\varepsilon\in\U$
and $\alpha\in\ida_{I(\mu)}^{-2}$.

\section{The Fundamental Domain\label{sec:The-Fundamental-Domain}}

The fundamental domain we describe in this section is essentially
the same as that used by Bluhmenthal \cite{MR1511187}, Mass \cite{MR0003405},
Tamagawa \cite{MR0121356}, Siegel \cite{MR659851} and others. The
main difference in these authors' approaches is in the description
of the ``bottom'' part which consists of a collection of hypersurfaces.
Here we adopt the description given by Siegel \cite{MR659851} since
it is easy to use for the explicit reduction algorithm. We have aimed
to provide sufficient details to demonstrate the appropriateness and
correctness of the algorithm and refer to \cite{MR659851} for details
and proofs. 

\subsection{Reduction with respect to units\label{subsec:Reduction-units}}

We use $\log:\mathbb{R}^{+}\rightarrow\mathbb{R}$ to denote the
natural logarithm and without risk of confusion we use the same notation
for the extended map $\log:\RK^{+}\rightarrow\mathbb{R}^{n}$ defined
by $\log(\vec{x})=\left(\log\emb_{1}(\vec{x}),\ldots,\log\emb_{n}(\vec{x})\right)$.
It is immediate from Dirichlet's unit theorem that the group of units
squared, $\U^{2},$ corresponds to an integral lattice $\Lambda$
of rank $n-1$ in $\mathbb{R}^{n}$, explicitly given by: 
\[
\Lambda=\log(\U^{2})=\left\{ \log(\varepsilon)\,:\,\varepsilon\in\U^{2}\right\} =\left\{ \sum_{k=1}^{n-1}a_{k}\log(|\varepsilon_{k}|)\,:\:a_{k}\in2\mathbb{Z}\right\} .
\]
The vectors $\log|\varepsilon_{k}|=(\log|\emb_{1}\varepsilon_{k}|,\ldots,\log|\emb_{n}\varepsilon_{k}|)^{t}$
form a basis of $\Lambda$ and we let \textbf{$B_{\Lambda}=(b_{rk})_{1\le r\le n,1\le k\le n-1}\in M_{n\times n-1}(\mathbb{R})$
}with $b_{rk}=\log\left|\emb_{r}(\varepsilon_{k})\right|$ denote
the corresponding basis matrix. Since all units have norm $1$ it
is easy to see that $\Lambda$ is contained in the $n-1$-dimensional
hyperplane 
\[
H=\left\{ \vec{u}\in\mathbb{R}^{n}\mid u_{1}+\ldots+u_{n}=0\right\} .
\]
We follow the explicit construction by Siegel and choose $K_{\Lambda}=B_{\Lambda}[-1,1[^{n-1}$
as a fundamental parallelepiped for $\Lambda$ and say that a vector
in $H$ is $\Lambda$-reduced if it belongs to $K_{\Lambda}$. If
$\vec{y}\in\RK^{+}$ we define $\tilde{\vec{y}}=\vec{y}\cdot(\N\vec{y})^{-1/n}$
and observe that $\N\tilde{\vec{y}}=1$, hence $\log(\tilde{\vec{y}})\in H$
and we say that $\vec{y}$ is $\U^{2}$-reduced if $\log(\tilde{\vec{y}})$
is $\Lambda$-reduced. This means that we can write 
\begin{equation}
B_{\Lambda}\vec{Y}=\log(\tilde{\vec{y}}),\label{eq:Y-reduced}
\end{equation}
where $\vec{Y}\in[-1,1[^{n-1}$ . The complete coordinate map $\vec{Y}_{\Lambda}:\mathbb{R}_{K}^{+}\rightarrow[-1,1[^{n-1}$
is then defined by setting $\vec{Y}_{\Lambda}(\vec{y})=\vec{Y}$ where
$\vec{Y}$ satisfies \eqref{eq:Y-reduced}. Observe that if $\vec{b}=(b_{1},\ldots,b_{n-1})^{t}\in\mathbb{Z}^{n-1}$
then 
\[
B_{\Lambda}\vec{b}=\sum_{k=1}^{n-1}b_{k}\log|\varepsilon_{k}|
\]
and therefore, if $\varepsilon=\varepsilon_{1}^{b_{1}}\cdots\varepsilon_{n-1}^{b_{n-1}}\in\U$
then $\vec{Y}_{\Lambda}(\varepsilon^{2}\vec{y})=\vec{Y}_{\Lambda}(\vec{y})+2\vec{b}.$
It follows that if we are given a $\vec{y}\in\mathbb{R}_{K}^{+}$
with $\vec{Y}_{\Lambda}(\vec{y})=(Y_{1},\ldots,Y_{n-1})$ and choose
$b_{i}=-\left\lfloor \frac{Y_{i}}{2}+\frac{1}{2}\right\rfloor $ then
$E(\varepsilon)\vec{y}=\varepsilon^{2}\vec{y}$ will be $\U^{2}$
reduced. Here $\left\lfloor x\right\rfloor $ is the nearest integer
to $x$, defined as the unique integer $n$ satisfying $x-1/2<n\le x+1/2$. 

It is easy to see that $K_{\Lambda}\simeq\Lambda\backslash\mathbb{R}^{n-1}$
is isomorphic via the logarithm map to a fundamental domain for the
action of the set of multipliers $\M_{\U}$ on $\RK^{+}$. For the
explicit computations of reduced vectors it is useful to have the
following explicit estimates in terms of the absolute row sums of
$B_{\Lambda}$: 
\[
r_{i}(B_{\Lambda})=\sum_{j=1}^{n-1}\left|\log\left|\emb_{i}\varepsilon_{j}\right|\right|,\quad1\le i\le n,
\]
and we observe that $\left\Vert B_{\Lambda}\right\Vert _{\infty}=\max s_{i}(B_{\Lambda}).$
\begin{lem}
If $\vec{u}\in\mathbb{R}^{n}$ is $\Lambda$-reduced then $|u_{i}|\le r_{i}(B_{\Lambda})$.
\end{lem}

\begin{proof}
If $\vec{u}\in K_{\Lambda}$ then $\vec{u}=B_{\Lambda}\vec{Y}$ for
some $\vec{Y}\in[-1,1]^{n-1}$ and hence 
\[
|u_{i}|=(B_{\Lambda}\vec{Y})_{i}\le\sum_{j=1}^{n-1}\left|\log\left|\emb_{i}\varepsilon_{j}\right|\right|\left|Y_{j}\right|\le r_{i}(B_{\Lambda}).
\]
\end{proof}
The following corollary is now immediate. 
\begin{cor}
\label{cor:mult_by_unit_reduced}If $\vec{y}\in\RK^{+}$ then there
is a unit $\varepsilon\in\U^{2}$ such that 
\[
(\N\vec{y})^{1/n}e^{-r_{i}(B_{\Lambda})}\le\left|\emb_{i}(\varepsilon\vec{y})\right|\le(\N\vec{y})^{1/n}e^{r_{i}(B_{\Lambda})}\quad\text{for all }1\le i\le n.
\]
\end{cor}

\subsection{Reduction with respect to translations\label{subsec:Reduction-translations}}

Let $\ida$ be an integral ideal in $\OK$ and choose an integral
basis $\beta_{1}^{(\ida)},\ldots,\beta_{n}^{(\ida)}$ of $\ida$.
Using the embedding map we identify $\ida$ with a lattice of rank
$n$ in $\mathbb{R}^{n}$, also denoted by $\ida$. The basis matrix
for this lattice is denoted by $B_{\ida}$ and we choose a fundamental
polytope $K_{\ida}=B_{\ida}[-1/2,1/2[^{n}$. For an element $\vec{x}\in\RK$
we define the $\ida$-coordinate vector $\vec{X}_{\ida}(\vec{x})$
by the equation 
\[
B_{\ida}\vec{X}_{\ida}(\vec{x})=\emb(\vec{x})
\]
and say that $\vec{x}$ is $\ida$-reduced if $\emb(\vec{x})\in K_{\ida},$
or in other words, if $\vec{X}_{\ida}(\vec{x})=(X_{1},\ldots,X_{n})$
with $-1/2\le X_{k}<1/2$ for all $k$s. 

If $\alpha=\sum_{k=1}^{n}a_{k}\beta_{k}^{(\ida)}\in\ida$ then $\vec{X}_{\ida}(\alpha)=(a_{1},\ldots,a_{n})$
and it is clear that $\vec{X}_{\ida}(\vec{x}+\alpha)=(X_{1}+a_{1},\ldots,X_{n}+a_{n})$
and hence, if we choose $a_{k}=-\left\lfloor X_{k}\right\rfloor $
then $T^{\alpha}\vec{x}=\vec{x}+\alpha$ will be $\ida$-reduced.

\subsection{Fundamental domain for the cusp stabilizer}

Let $\lambda\in\PK$ be a cusp of $\Gamma_{K}$, $\ida_{\lambda}$
the corresponding representative ideal and $\idb=\ida_{\lambda}^{-2}$
with an integral basis $\beta_{1}^{(\idb)},\ldots,\beta_{n}^{(\idb)}.$
For an element $\vec{z}\in\HK$ we define $\vec{z}_{\lambda}=\vec{x}_{\lambda}+i\vec{y}_{\lambda}=A_{\lambda}^{-1}\vec{z}$
and say that $\vec{z}$ is reduced with respect to $\lambda$ if $\vec{x}_{\lambda}$
is reduced with respect to $\idb$ and $\vec{y}_{\lambda}$ is reduced
with respect to $\U^{2}$. We let $\G_{\lambda}$ denote the set of
all such reduced points, more precisely 
\[
\G_{\lambda}=\left\{ \vec{z}\in\HK\mid\vec{X}_{\idb}(\vec{x}_{\lambda})\in[-1/2,1/2[^{n}\,\text{and}\ \vec{Y}_{\Lambda}(\vec{y}_{\lambda})\in[-1/2,1/2[^{n-1}\right\} .
\]
It is easy to show that the set $\G_{\lambda}$ is indeed a fundamental
domain for the action of $\Gamma_{K,\lambda}=A_{\lambda}\Gamma(\OK\oplus\ida_{\lambda}^{2})A_{\lambda}^{-1}$
on $\HK$. Note that for the modular group, $\PSLZ$, the analogue
of the domain $\G_{\lambda}$ is the strip $-1/2<\Re(z)\le1/2$.

\subsection{Cuspidal regions}

If the regions $\G_{\lambda}$ in the previous section are analogues
of the vertical strip we will now look at the analog of the curved
part of the fundamental domain, given by $|z|\ge1$. For the modular
group this can be interpreted in terms of a reflection in the isometric
circle corresponding to the map given by $z\mapsto-z^{-1}$. An analog
interpretation is valid for Hilbert modular groups but it is much
harder to work out precisely which reflections to include even for
small number fields of class number $1$. 

If $\vec{z}=\vec{x}+i\vec{y}\in\HK$ we define $\Delta(\vec{z},\infty)$,
the \emph{distance} to the cusp at $\infty$, by
\[
\Delta(\vec{z},\infty)=\N\left(\Im\vec{z}\right)^{-1/2}
\]
and the distance to an arbitrary cusp $\mu=\cusp{\rho}{\sigma}$ with
associated ideal $\ida=(\rho,\sigma)$ is 
\begin{align}
\Delta(\vec{z},\mu) & =\N\left(\ida\right)^{-1}\N(\Im A_{\mu}^{-1}\vec{z})^{-1/2}=\N\left(\ida\right)^{-1}\N\left(\vec{y}\right)^{-1/2}\N(\left|-\sigma\vec{z}+\rho\right|^{2})^{1/2}\label{eq:dist_formula}\\
 & =\N(\ida)^{-1}\N\left(\left(-\sigma\vec{x}+\rho\right)^{2}\vec{y}^{-1}+\sigma^{2}\vec{y}\right)^{\frac{1}{2}},\nonumber 
\end{align}
where $\N(\ida)$ is the norm of the ideal $\ida$. This expression
is independent of the choice of representatives $\rho$ and $\sigma$
as well as the choice of $A_{\mu}$. Observe that the normalization
factor $N(\mathfrak{a})^{-1}$, which accounts for the independence
of the choices of $\rho$ and $\sigma$ is present in \cite{MR930101}
but not in \cite{MR659851}. The expression $\Delta(\vec{z},\mu)$
is in fact bi-invariant under $\Gamma_{K}$, in other words, $\Delta\left(A\vec{z},A\mu\right)=\Delta\left(\vec{z},\mu\right)$
for all $A\in\Gamma_{K}$. We will show later, in Lemma \ref{lem:closest_cusp_exists},
that for every $\vec{z}\in\HK$ there exists a cusp $\lambda$ which
is closest to $\vec{z}$ and it follows that the\emph{ invariant height}

\[
\Delta(\vec{z})=\inf\left\{ \Delta(\vec{z},\lambda)\mid\lambda\in\PK\right\} 
\]
is well-defined, invariant under $\Gamma_{K}$ and $\Delta(\vec{z})=\Delta(\vec{z},\lambda)$
for some cusp $\lambda$ (not necessarily unique).

We are now fully prepared to give the definition of the fundamental
domain that we are interested in. For a cusp representative $\lambda_{j}$
with $1\le j\le h$ we let \emph{$\mathcal{F}_{j}$ }denote the set
of $\lambda_{j}$-reduced points that are closest to $\lambda_{j}$,
in other words:
\[
\mathcal{F}_{j}=\left\{ \vec{z}\in\G_{\lambda_{j}}\mid\Delta\left(\vec{z}\right)=\Delta(\vec{z},\lambda_{j})\right\} .
\]
It can then be shown (cf. e.g. \cite{MR659851}) that the set 
\[
\mathcal{F}_{K}=\cup_{j=1}^{h}\mathcal{F}_{j}
\]
is a fundamental domain for the action of $\Gamma_{K}$ on $\HK$.
Given that $\mathcal{F}_{K}$ is a fundamental domain we now turn
to the problem of reducing a point $\vec{z}$ to its representative
inside $\mathcal{F}_{K}.$ It is clear that the as soon as we find
a closest cusp, say $\mu\in\PK$, which is equivalent to a cusp representative
$\lambda_{j}$ with $U_{\mu}(\mu)=\lambda_{j}$ then $\lambda_{j}$
is a closest cusp to $U_{\mu}\vec{z}$ and we can use the straight-forward
reduction with respect to units and translations from Sections \ref{subsec:Reduction-units}
and \ref{subsec:Reduction-translations} to find an $\varepsilon\in\U^{2}$
and $\alpha\in\ida_{j}^{-2}$ such that $\vec{z}^{*}=A_{\lambda_{j}}T^{\alpha}E(\varepsilon)A_{\lambda_{j}}^{-1}U_{\mu}\vec{z}$
belongs to $\mathcal{F}_{j}$. 

The reduction with respect to units and translations is essentially
done in constant time independent of $\vec{z}$ and has been efficiently
implemented by many authors, cf. e.g. \cite{MR3376741}. The practical
and theoretical complexity of the reduction algorithm is almost entirely
in the finding of the closest cusp. The next section is dedicated
to auxiliary results and details on how our algorithm for finding
the closest cusp works and we will then summarize the actual algorithm
in the following section.

\section{Finding the closest cusp\label{sec:Finding-the-closest}}

Our approach to finding the closest cusp $\lambda$ is to analyze
the existence and conditional uniqueness proofs from the lecture notes
of Siegel \cite{MR659851} and find explicit and efficient bounds
for all constants involved. The general idea was already present in
a slightly different form in the work of Maaß\cite{MR0003405} but
note that some of the explicit constants present, in e.g. Hilfssatz
II, are in general weaker than those we obtain here. The aim of this
section is to include sufficient details in the proofs for a reader
to be able to both understand and verify the functionality of the
associated code \cite{hilbertmodgrp} as well as being able to implement
these algorithms independently. 
\begin{lem}
\label{lem:closest_cusp_exists}If $\vec{z}\in\HK$ then there exists
a cusp $\lambda\in\PK$ such that 
\[
\Delta(\vec{z},\lambda)\le\Delta(\vec{z},\mu)\quad\forall\mu\in\PK.
\]
\end{lem}

\begin{proof}
Let $\vec{z}=\vec{x}+i\vec{y}\in\HK$ be fixed. It is sufficient to
show that for any given cusp $\mu$ there exists only a finite number
of cusps $\lambda$ such that $\Delta(\vec{z},\lambda)\le\Delta(\vec{z},\mu)$. 

It follows from Section \ref{subsec:Cusps-of-Hilbert} that we can
assume that $\lambda=\cusp{\rho}{\sigma}$ where $\rho$ and $\sigma$
are chosen such that $(\rho,\sigma)=\ida_{i}$ for some class group
representative $\ida_{i}$ and in particular $\N((\rho,\sigma))\le C$
where 
\[
C=\max\left\{ \N(\ida_{1}),\ldots,\N(\ida_{h})\right\} .
\]
We now consider $\Delta(\vec{z},\lambda)$ as a function of the algebraic
integers $\sigma$ and $\rho$ and write
\[
\Delta_{\vec{z}}(\rho,\sigma)\coloneqq\Delta(\vec{z},(\rho:\sigma))=\N((\rho,\sigma))^{-1}\left(\N\left(\left(-\sigma\vec{x}+\rho\right)^{2}\vec{y}^{-1}+\sigma^{2}\vec{y}\right)\right)^{\frac{1}{2}}.
\]
It is sufficient to show that if $d>0$ there exists only a finite
number of pairs $\rho,\sigma\in\OK$ modulo units, such that $\Delta_{\vec{z}}\left(\rho,\sigma\right)<d$.
Given such a pair write
\[
\Delta_{\vec{z}}(\rho,\sigma)=\N((\rho,\sigma))^{-1}(\N\vec{w})^{1/2},
\]
where $\vec{w}=\left(-\sigma\vec{x}+\rho\right)^{2}\vec{y}^{-1}+\sigma^{2}\vec{y}\in\mathbb{R}_{K}^{+}$.
It follows from Corollary \ref{cor:mult_by_unit_reduced} that there
exists a unit $\varepsilon\in\U$ such that 
\[
\left|\emb_{i}\left(\varepsilon^{2}\vec{w}\right)\right|\le e^{r_{i}(B_{\Lambda})}(\N\vec{w})^{\frac{1}{n}}\le e^{r_{i}(B_{\Lambda})}d^{2/n}C^{2/n},\quad\text{for all}\quad1\le i\le n.
\]
Setting $\delta_{i}=e^{r_{i}(B_{\Lambda})}d^{2/n}C^{2/n}$ we can
therefore assume that $\sigma$ and $\rho$ have been chosen such
that $|\emb_{i}\left(\vec{w}\right)|\le\delta_{i}$, and hence that
\[
\emb_{i}(\sigma^{2}\vec{y})\le\delta_{i}\quad\text{and}\quad\emb_{i}(\left(-\sigma\vec{x}+\rho\right)^{2}\vec{y}^{-1})\le\delta_{i}.
\]
It follows that the coordinates of the embeddings of $\sigma$ and
$\rho$ are bounded by 
\begin{align}
\left|\sigma_{i}\right|^{2} & \le\delta_{i}y_{i}^{-1}\quad\text{and}\label{eq:ineq-sigma}\\
\left|\rho_{i}-\sigma_{i}x_{i}\right|^{2} & \le\delta_{i}y_{i}.\label{eq:ineq-rho}
\end{align}
The inequalities \eqref{eq:ineq-sigma} and \eqref{eq:ineq-rho} clearly
define a bounded domain in $\mathbb{R}^{n}\times\mathbb{R}^{n}$ and
the statement follows since the embeddings of $\OK$ form a lattice
in $\mathbb{R}^{n}$.
\end{proof}
An immediate consequence of the previous proof, and in particular
\eqref{eq:ineq-sigma} and \eqref{eq:ineq-rho} is the following result
which is crucial to our algorithm. 
\begin{thm}
\label{thm:sigma_rho_bound_emb}Let $\vec{z}\in\HK$ and assume that
there is a cusp $\lambda$ with $\Delta(z,\lambda)=d$. Then a closest
cusp can be chosen as $(\rho:\sigma)$ where the embeddings of $\rho$
and $\sigma$ satisfy the following bounds: 
\end{thm}

\[
|\sigma_{i}|\le D_{i}\cdot d^{1/n}y_{i}^{-1/2}\quad\text{ and}\quad\left|\rho_{i}-x_{i}\sigma_{i}\right|\le D_{i}\cdot d^{1/n}y_{i}^{1/2},
\]
where 
\[
D_{i}=C^{1/n}e^{\frac{1}{2}r_{i}(B_{\Lambda})},
\]
and, additionally, the norms are bounded by 
\[
\N\left(\left|\sigma\right|\right)\le dC\N\left(\vec{y}\right)^{-1/2}\quad\text{and}\quad\N\left(\left|-\sigma\vec{x}+\rho\right|\right)\le dC\N\left(\vec{y}\right)^{1/2}.
\]

To apply the previous theorem we need to find an initial cusp $\lambda$.
It is, for instance, always possible to choose $\infty$, in which
case $d=\Delta\left(z,\infty\right)=\N(\vec{y})^{-1/2}$, or $0$,
in which case $d=\Delta\left(z,0\right)=\N(\vec{y})^{-1/2}\N(\vec{x}^{2}+\vec{y}^{2})^{1/2}$
. However, it is clear that we would like to obtain as small initial
bound as possible and if $\N(\vec{y})$ is small than we need to find
another cusp to start with. 

Fortunately there is a method which seems to work well in practice
when $\N\left(\vec{y}\right)$ is small. This method was introduced
by Bouyer and Streng \cite{MR3376741} and the main idea is to use
LLL reduction to find a vector of short norm, $-\sigma\vec{z}+\rho$,
in the lattice $L_{\vec{z}}=\OK\vec{z}+\OK$ and the corresponding
cusp $(\rho:\sigma)$ will then be close to $\vec{z}$ by \eqref{eq:dist_formula}. 
\begin{rem}
It should be noted that the LLL reduction method by itself does not
does not necessarily yield the closest cusp, as the LLL algorithm
is not guaranteed to return the shortest vector and the definition
of distance $\Delta(z,(\rho:\sigma))$ also involves the norm of the
ideal $(\rho,\sigma)$. For a provably correct algorithm (in all degrees)
it is therefore necessary to combine this preliminary optimization
with an exhaustive search using the explicit bounds of Lemma \ref{thm:sigma_rho_bound_emb}.
\end{rem}

Since the only integer in $\OK$ with norm less than $1$ is $0$
the norm bound of Theorem \ref{thm:sigma_rho_bound_emb} immediately
implies the following. 
\begin{cor}
\label{cor:infinity_closest}Let $\vec{z}=\vec{x}+i\vec{y}\in\HK$.
If $\N(\vec{y})>C$ then $\infty$ is the closest cusp to $\vec{z}$.
\end{cor}

Unfortunately it is in general not so easy to find the closest cusp
and we will see that it is often necessary to compare distances to
many different cusps. However, the number of comparisons needed can
sometimes be reduced by using the following Lemma and Corollary. 
\begin{lem}
\label{lem:exists_d_cusp_same}There exists a constant $d>0$, depending
only on $K$, such that for all $\vec{z}=\vec{x}+i\vec{y}\in\HK$,
if $\lambda$ and $\mu$ are cusps of $K$ with $\Delta(\lambda,\vec{z})<d$
and $\Delta(\mu,\vec{z})<d$ then $\lambda=\mu$. 
\end{lem}

\begin{proof}
Let $\vec{z}=\vec{x}+i\vec{y}\in\HK$. Assume that $\lambda=(\sigma:\rho)$
and $\mu=(\sigma_{1}:\rho_{1})$ satisfy $\Delta(\lambda,\vec{z})<d$
and $\Delta(\mu,\vec{z})<d$ for some positive $d$. Observe that
the algebraic integer $\rho\sigma_{1}-\sigma\rho_{1}$ can be written
\begin{align*}
\rho\sigma_{1}-\sigma\rho_{1} & =\left(-\sigma\vec{x}+\rho\right)\vec{y}^{-1/2}\sigma_{1}\vec{y}^{1/2}-\left(-\sigma_{1}\vec{x}+\rho_{1}\right)\vec{y}^{-1/2}\sigma\vec{y}^{1/2}.
\end{align*}
Since \eqref{eq:ineq-sigma} and \eqref{eq:ineq-rho} applies to both
$(\sigma,\rho)$ and $(\sigma_{1},\rho_{1})$ it is easy to see that
\[
\emb_{i}\left(\left|\rho\sigma_{1}-\sigma\rho_{1}\right|\right)\le2\delta_{i},
\]
where $\delta_{i}=e^{r_{i}(B_{\Lambda})}d^{2/n}C^{2/n}$. It follows
that $\N\left(\rho\sigma_{1}-\sigma\rho_{1}\right)\le2^{n}\prod\delta_{i}$
and hence, if $d<C^{-1}2^{-n/2}e^{-\sum r_{i}(B_{\Lambda})}$ then
we must have $\rho\sigma_{1}-\sigma\rho_{1}=0$ so $\mu=\lambda$. 
\end{proof}
\begin{cor}
\label{cor:estimate_closest_cusp}Let $\vec{z}=\vec{x}+i\vec{y}\in\HK$.
If $\lambda$ is a cusp with $\Delta(\lambda,\vec{z})<C^{-1}2^{-n/2}e^{-\sum r_{i}(B_{\Lambda})}$
then $\lambda$ is the closest cusp to $\vec{z}$. 
\end{cor}

The previous lemma also has the geometric consequence that it is possible
to decompose the fundamental domain $\mathcal{F}$ into a compact
part and disjoint cuspidal parts.

\section{Algorithms\label{sec:Algorithms}}

We will now describe the actual reduction algorithm in detail. The
key idea is to use Theorem \ref{thm:sigma_rho_bound_emb} to find
bounded regions in $\mathbb{R}^{n}$ where the embeddings of the numerators
and denominators of potential closest cusps must be located. We then
compare the distance to $\vec{z}$ for each of the candidate cusps,
except if one of the distances is less than the bound in Corollaries
\ref{cor:infinity_closest} or \ref{cor:estimate_closest_cusp}, in
which case we terminate the search early. 

Recall that we have a fixed integral basis $\alpha_{1},\ldots,\alpha_{n}$
of $\OK$ and a corresponding lattice in $\mathbb{R}^{n}$ with basis
matrix $B_{\OK}$. If $\beta\in\OK$ is given by $\beta=\sum_{i=1}^{n}X_{i}\alpha_{i}$
for some integer vector $\vec{X}\in\mathbb{Z}^{n}$ then the embeddings
of $\beta$ correspond to the vector $\emb(\beta)=B_{\OK}\vec{X}$
in $\mathbb{R}^{n}$. If we can bound the vector $\emb$($\beta)$
in a parallelotope $P$ it follows that $\vec{X}$ must belong to
the polytope $B_{\OK}^{-1}(P)$ and we thus need to find vectors with
integer coordinates inside this set. 

A preliminary investigation of the performance showed that the most
efficient way to find these seems to be to search for integer vectors
in a bounding parallelotope of $B_{\OK}^{-1}(P)$, which we denote
by $\mathbb{BP}(B_{\OK}^{-1}(P))$, and apply the embedding map to
test whether or not to include them in the result. Using this idea
together with Theorem \ref{thm:sigma_rho_bound_emb} gives us the
following algorithm.
\begin{algorithm}
[Finding the closest cusp]\label{alg:closest-cusp-2}Let $K$ be
a fixed totally real number field, all notation be as above and let
$\vec{z}=\vec{x}+i\vec{y}\in\HK.$
\begin{lyxlist}{00.00.0000}
\item [{Step~1:}] If $\N(\vec{y})>C$ return $\infty=(0:1)$ as the closest
cusp. 
\item [{Step~2:}] Use the LLL reduction to find a potentially closest
cusp, $\lambda$, and set $d=\min\left\{ \Delta(\vec{z},\lambda),\Delta(\vec{z},\infty),\Delta(\vec{z},0)\right\} $ 
\item [{Step~3:}] Recall that $C=\max\left(\N(\mathfrak{a}_{1}),\ldots\N(\mathfrak{a}_{h})\right)$
and $D_{i}=\max C^{1/n}e^{\frac{1}{2}r_{i}(B_{\Lambda})}$. Define
\[
a_{i}=D_{i}d^{1/n}y_{i}^{-1/2},\ 1\le i\le n,\quad\text{and set}
\]
\[
P_{\sigma}=[-a_{1},a_{1}]\times\cdots\times[-a_{n},a_{n}]\quad\text{and}\quad\hat{P}_{\sigma}=\mathbb{B}\mathbb{P}(B_{\OK}^{-1}(P_{\sigma})).
\]
\item [{Step~4:}] Compute the integral points $\vec{X}^{(1)},\ldots,\vec{X}^{(M)}$
of $\hat{P}_{\sigma}$, and for each $1\le j\le M$: 
\begin{lyxlist}{00.00.0000}
\item [{(a)}] Compute the corresponding $\emb(\sigma)=B_{\OK}\vec{X}_{j}$
and if 
\[\N\left(\sigma\right)>dC\N\left(\vec{y}\right)^{-1/2} \quad \textrm{or}\quad  \emb(\sigma)\notin P_{\sigma},\]
 remove the corresponding $\vec{X}^{(j)}$
from the list and repeat for the next $j$, if not, go to the next
step.
\item [{(b)}] Set $b_{j,i}^{\pm}=\sigma_{i}x_{i}\pm y_{i}a_{i}$, and let
\[
P_{\rho,j}=[b_{j,1}^{-},b_{j,1}^{+}]\times\cdots\times[b_{j,n}^{-},b_{j,n}^{+}]\quad\text{and}\quad\hat{P}_{\rho,j}=\mathbb{B}\mathbb{P}(B_{\OK}^{-1}(P_{\rho,j})).
\]
\item [{(c)}] Compute the integral points $\vec{Y}^{(j,1)},\ldots,\vec{Y}^{(j,N(j))}$
of $\hat{P}_{\rho,j},$ and for each $1\le i\le N(j),$ compute $\emb(\rho)=B_{\OK}\vec{Y}^{(j,i)}$
and if $\emb(\rho)\notin P_{\rho,j}$ remove the corresponding $\vec{Y}^{(j,i)}$
from the list. 
\end{lyxlist}
\end{lyxlist}
After relabeling the remaining vectors if necessary we find that a closest
cusp to $\vec{z}$ can now be found corresponding to a pair in the
finite set
\[
\left\{ \left(\rho,\sigma\right)\mid\sigma=B_{\OK}\vec{X}^{(j)},\rho=B_{\OK}\vec{Y}^{(j,k)},\ 1\le j\le M',\,1\le k\le N'(j)\right\} 
\]
where $M'$ and $N'(j)$ are some positive integers. 
\end{algorithm}

\begin{rem}
Note that we do not make explicit use of the norm bound for $\rho$
here, it is instead part of finding the minimal distance in the final
set. 
\end{rem}

We can now combine Algorithm \ref{alg:closest-cusp-2} with the reduction by units and translation described in Section \ref{sec:The-Fundamental-Domain} 
to formulate the complete reduction algorithm. 
\begin{algorithm}
\label{alg:Reduction-algorithm}[Reduction algorithm]Let $K$ be
a fixed totally real number field, let $\vec{z}\in\HK$ and assume
all notation is as above, 
\begin{lyxlist}{00.00.0000}
\item [{Step~1:}] Use Algorithm \ref{alg:closest-cusp-2} to find the
closest cusp to $\vec{z}$, say $\mu$.
\item [{Step~2:}] Find the cusp representative, $\lambda_{j}$, corresponding
to $\mu$ and $U_{\mu}\in\Gamma_{K}$ such that $U_{\mu}(\mu)=\lambda_{j}$
. 
\item [{Step~3:}] Set $\vec{z}_{\lambda_{j}}=A_{j}^{-1}U_{\mu}\vec{z}=\vec{x}_{\lambda_{j}}+i\vec{y}_{\lambda_{j}}$. 
\item [{Step~4:}] Let $\vec{Y}=\vec{Y}_{\Lambda}(\vec{y}_{\lambda_{j}})$
and define $\varepsilon=\varepsilon_{1}^{b_{1}}\ldots\varepsilon_{n-1}^{b_{n-1}}$
where $b_{k}=-\left\lfloor \frac{Y_{k}}{2}\right\rfloor $. 
\item [{Step~5:}] Set $\vec{z}'=E(\varepsilon)\vec{z}_{\lambda_{j}}=\vec{x}'+i\vec{y}'$.
\item [{Step~6:}] Let $\vec{X}=\vec{X}_{\ida_{j}^{-2}}(\vec{x}')$ and
define $\alpha=a_{1}\beta_{1}+\cdots+a_{n}\beta_{n}$ where $a_{k}=-\left\lfloor X_{k}\right\rfloor $
and $\beta_{1},\ldots,\beta_{n}$ is an integral basis for $\ida_{j}^{-2}$. 
\item [{Step~7:}] Set $A=A_{j}T^{\alpha}E(\varepsilon)A_{j}^{-1}U_{\mu}$. 
\end{lyxlist}
\end{algorithm}

Then $A\in\Gamma_{K}$ and $A\vec{z}\in\mathcal{F}_{j}\subseteq\mathcal{F}$.

\subsection{A brief analysis of runtime and performance}

 It is clear that reduction within the cuspidal domain is essentially
of constant time with respect to $\vec{z}$. The run-time is therefore
essentially proportional to the total number of potential $\sigma$s
and $\rho$s that are investigated in Algorithm \eqref{alg:closest-cusp-2}
and each of these numbers are proportional to the volumes of the corresponding
polytopes. It is of little practical use to make a very precise run-time
analysis here but by using appropriate upper bounds it is easy to
see that for a fixed totally real number field $K$ the run-time is
polynomial in $\left\Vert \vec{x}\right\Vert _{\infty},\left\Vert \vec{y}\right\Vert _{\infty}$
and $\left\Vert \smash{\vec{y}^{-1}}\right\Vert _{\infty}$ as $\vec{z}$
varies. Similarly, if $\vec{z}$ is fixed and we let $K$ vary then
the run-time is exponential in the degree of $K$ and $\left\Vert B_{\Lambda}\right\Vert _{\infty}$,
and polynomial in $C$, $\left\Vert B_{\OK}^{-1}\right\Vert _{\infty}$
and $\left\Vert B_{\OK}\right\Vert _{\infty}$. While algebraic quantities
like the discriminant and regulator of $K$ do play a direct role also in
the reduction by units and translations, these can be bounded by the
respective matrix norms. 

While a more precise analysis for the dependency on $\vec{z}$ is
not too difficult to perform, a detailed analysis on the precise dependency
on the number field is more complex due to the number of different
parameters involved. For testing the runtime in practice it is convenient
to consider the point $\vec{z}=i\vec{1}$ since it will always be
closest to both $0$ and $\infty$ and the preliminary search using
LLL does not provide any better bound. Table \ref{tab:Times-exa}
contains times to find the closest cusp of $i\vec{1}$ for the different
fields we consider in Section \ref{sec:Examples}. Here $\alpha_{1}$
and $\alpha_{2}$ have minimal polynomials $\alpha_{1}^{3}-\alpha_{1}^{2}-2\alpha_{1}+1$
and $\alpha_{2}^{3}-\alpha_{2}^{2}-2\alpha_{2},$ respectively. 
For a more systematic comparison regarding the dependency on the discriminant
we also compared quadratic fields of class number one and discriminant up to $100$. See Table \ref{tab:Times-exa2}. 
The difference in timing between discriminant $93$ and $97$ of a factor over $200$ is striking. 
It highlights that the influence of the discriminant is vastly overshadowed by that of the size of the embeddings  of the fundamental units. 
The lengths in question here are $\approx 3.37$ and $\approx 9.32$,  respectively, and $\exp( 9.3 - 3.3)\approx	403$.

All computations below were performed on a single 2GHz Xeon E5-2660
core and the reported time is an average of 100 runs. 
\begin{table}
\caption{\label{tab:Times-exa}Times to find closest cusp of $\vec{z}=i\vec{1}$ for different fields.}

\begin{center}

\begin{tabular}{|c|c|c|c|c|}
\hline 
Field & $\mathbb{Q}(\sqrt{5})$ & $\mathbb{Q}(\sqrt{10})$ & $\mathbb{Q}(\alpha_{1})$ & $\mathbb{Q}(\alpha_{2})$ \tabularnewline
\hline 
\hline 
Time / ms& 13ms & 19ms & 124ms & 387s\tabularnewline
\hline 
\end{tabular}
\end{center}

\end{table}
 \begin{table}
\caption{\label{tab:Times-exa2}Times to find closest cusp of $\vec{z}=i\vec{1}$ for quadratic fields of discriminant $D$ and class number $1$}
\begin{tabular}{|c|c|c|c|c|c|c|c|c|c|c|c|c|c|c|c|c|c|c|c|c|c|c|c|c|c|c|}
\hline 

D& 5 & 8 & 12 & 13 & 17 & 21 & 24 & 28 & 29 & 33 & 37 & 41 & 44  \\
\hline 
Time / ms& 12 & 14 & 14 & 13 & 17 & 14 & 18 & 20 & 14 & 37 & 17 & 47 & 23 \\
\hline 
\hline 
D& 53 & 56 & 57 & 61 & 69 & 73 & 76 & 77 & 88 & 89 & 92 & 93 & 97\\
\hline 
Time / ms&15 & 26 & 167 & 30 & 23 & 1087 & 178 & 16 & 201 & 488 & 33 & 24 & 4740\\
\hline 
\end{tabular}

\end{table}

\subsection{A note on the implementation}

The algorithms described above are currently implemented as part of
a package in Python with parts written in Cython and is dependent
on SageMath \cite{sage}. The package is available from \cite{hilbertmodgrp}
and open sourced under GPLv3+. 

\section{\label{sec:Examples}Examples}

The aim of the examples presented here is to demonstrate how the algorithm
works as well as making it easy for readers to verify the correctness.
We will consider three examples in detail: first the standard example
of $K_{1}=\mathbb{Q}(\sqrt{5})$ , which has degree $2$, discriminant
$5$ and class number $1$, then $K_{2}=\mathbb{Q}(\sqrt{10})$ which
has degree $2$, discriminant $40$ and class number $2$, followed
by $K_{3}=\mathbb{Q}(\alpha)$ where $\alpha$ has minimal polynomial
$\alpha^{3}-\alpha^{2}-2\alpha+1$, which has degree $3$, discriminant
$49$ and class number $1$. The computations involved in these three
examples are demonstrated in the accompanying Jupyter notebooks that
can be found in \cite{hilbertmodgrp}
\begin{example}
Consider $K=\mathbb{Q}\left(\sqrt{5}\right)$ with fundamental unit
$\varepsilon=\frac{1}{2}(1+\sqrt{5})$, ring of integers $\OK=\mathbb{Z}\oplus\mathbb{Z}\varepsilon$
and class number $1$. Here
\begin{align*}
B_{\Lambda} & =\left(\begin{array}{cc}
\log(\frac{1}{2}(1+\sqrt{5})) & \log(\frac{1}{2}(\sqrt{5}-1))\end{array}\right),\\
B_{\OK} & =\left(\begin{array}{cc}
1 & \frac{1}{2}\left(1+\sqrt{5}\right)\\
1 & \frac{1}{2}\left(1-\sqrt{5}\right)
\end{array}\right),\ B_{\OK}^{-1}=\frac{1}{-\sqrt{5}}\left(\begin{array}{cc}
\frac{1}{2}(1-\sqrt{5}) & -\frac{1}{2}(1+\sqrt{5})\\
-1 & 1
\end{array}\right)
\end{align*}
and it is immediate to see that 
\begin{align*}
r_{1}(B_{\Lambda})=r_{2}(B_{\Lambda}) & \approx0.48,\quad D_{1}=D_{2}\approx1.27,\quad\left\Vert B_{\OK}\right\Vert _{\infty}\approx2.62,\quad\text{and }\quad\left\Vert B_{\OK}^{-1}\right\Vert _{\infty}=1.
\end{align*}
Consider now Algorithm \ref{alg:closest-cusp-2} applied to $\vec{z}=\vec{y}=i\vec{1}\in\Hn$.
Since $y_{1}=y_{2}=1$ the first bounds are given by $a_{1}=a_{2}=D_{1}$
and it can be computed that $P_{\sigma}$ is the polygon bounded by
the vertices 
\[
B_{\OK}^{-1}((\pm D_{0},\pm D_{0}))=\{(0.57,1.14),(1.27,0.0),(-0.57,-1.14),(-1.27,0.0)\}.
\]
For this $\vec{z}$ the preliminary reduction does not produce any
better cusp than $\infty$ and the norm bound is given by $C\N\left(\vec{y}\right)^{-1/2}=1$.
The domain $P_{\sigma}$ and its pre-image together with the embedded
points of $\OK$ and $\mathbb{Z}^{2}$ and the curves indicating the
norm bound are shown in Figure \ref{fig:domainsK1}. Note that we
show the actual domain $B_{\OK}^{-1}(P_{\sigma})$ and not the bounding
box for extra clarity. 

\begin{figure}
\begin{minipage}[t]{0.49\linewidth}
\centering\resizebox{\linewidth}{!}{\includegraphics[scale=0.25]{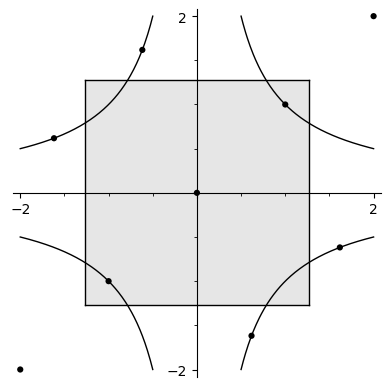}}
\end{minipage}     
\begin{minipage}[t]{0.49\linewidth}
\resizebox{\linewidth}{!}{\includegraphics[scale=0.25]{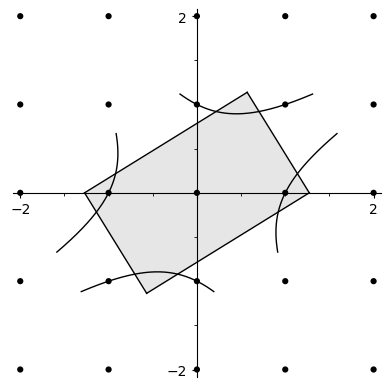}}
\end{minipage}
\caption{\label{fig:domainsK1}$B_{\mathcal{O}_K}(P_{\sigma}$)  and $P_{\sigma}$ together with $\N(\sigma)=1$ for $\mathbf{z}=i\mathbf{1}$ and $K=\mathbb{Q}(\sqrt{5})$}
\end{figure}

It is clear from the figure that the only integral points in $P_{\sigma}$
are $\vec{X}^{(1)}=(0,0)$, $\vec{X}^{(2)}=(1,0)$ and $\vec{X}^{(3)}=(-1,0)$
so the three candidates for $\sigma$ are $0$, $1$ and $-1$. The
value $\sigma=0$ of course corresponds to the cusp at infinity so
we can choose $\rho=1$ in this case. For $\sigma=\pm1$ we get $P_{\rho,2}=P_{\sigma}$
so the candidates for $\rho$ are $0$ and $\pm1$. The potential
cusps are therefore $\infty=(1:0)$, $c_{1}=(0:1)$, $c_{2}=(1:1)$
and $c_{3}=(-1:1)$. The corresponding distances are $\Delta(\vec{z},\infty)=\N(\vec{y})^{-1/2}=1,$
$\Delta(\vec{z},(0:1))=\N\left(\vec{y}\right)^{1/2}=1$ and 
\begin{align*}
\Delta(\vec{z},(1:1)) & =\N\left(\vec{y}^{-1}+\vec{y}\right)^{1/2}=2.
\end{align*}
Therefore both $\infty=(1:0)$ and $0=(0:1)$ are closest cups. 

If we consider instead $\vec{z}=\frac{1}{2}i\vec{1}$ then then a
preliminary search (using e.g. the LLL method) finds the cusp $0=(0:1)$
and it is easy to see that $\Delta(\vec{z},\infty)=2$ and $\Delta(\vec{z},0)=1/2$.
We can therefore apply the algorithm with an initial estimate of $d=1/2.$
This leads to the same bounds for $\sigma:$ $|\sigma_{i}|\le D_{0}$
but the bounds for $\rho$ gets scaled: $\left|\rho_{i}\right|\le D_{0}/2\approx0.636$.
The candidate cusps are therefore simply $(1:0)=\infty$ and $(0:1)$
with the closest cusp being $(0:1).$

If we had not performed the initial search and instead simply used
$d=\Delta\left(\vec{z},\infty\right)=2$ for the initial bound we
would have obtained $9$ candidates for sigma and in the end $9$
candidates for closest cusp .
\end{example}

\begin{figure}
\begin{minipage}[t]{0.49\linewidth}
\centering\resizebox{\linewidth}{!}{\includegraphics[scale=0.25]{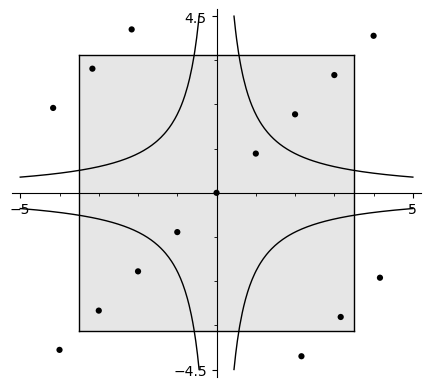}}

\end{minipage}     
\begin{minipage}[t]{0.49\linewidth}

\resizebox{\linewidth}{!}{\includegraphics[scale=0.25]{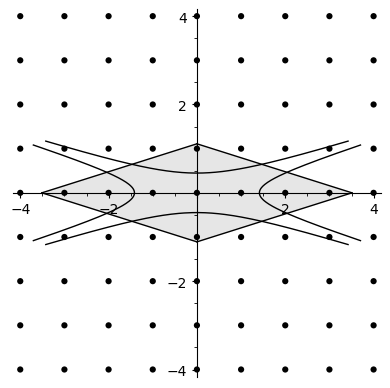}}
\end{minipage}

\caption{\label{fig:K2.z1.domains} $B_{\mathcal{O}^{-1}_{K_2}}(P_{\sigma})$  and $P_{\sigma}$ together with the curve $\N(\sigma)=2$ for $\mathbf{z}=i\mathbf{1}$, $K_2=\mathbb{Q}(\sqrt{10})$.}
\end{figure}

\begin{example}
Consider $K_{2}=\mathbb{Q}\left(\sqrt{10}\right)$ with fundamental
unit $\varepsilon=3+\sqrt{10}$, ring of integers $\OK=\mathbb{Z}\oplus\mathbb{Z}\sqrt{10}$
and class number $2$. The cusp representatives are 
\[
\lambda_{1}=\infty\quad\text{and}\quad\lambda_{2}=(2:\sqrt{10})
\]
 Note that the norm of the ideal associated with $\lambda_{2}$ is
$\N\left(\left(2,\sqrt{10}\right)\right)=2$. We now find that 
\begin{align*}
B_{\Lambda} & =\left(\begin{array}{cc}
\log((3+\sqrt{10})) & \log(\sqrt{10}-3))\end{array}\right),\\
B_{\OK} & =\left(\begin{array}{cc}
1 & \sqrt{10}\\
1 & -\sqrt{10}
\end{array}\right)\quad\text{and}\quad B_{\OK}^{-1}=\frac{1}{-2\sqrt{10}}\left(\begin{array}{cc}
-\sqrt{10} & -\sqrt{10}\\
-1 & 1
\end{array}\right),
\end{align*}
and it is immediate to see that 
\begin{align*}
r_{1}(B_{\Lambda})=r_{2}(B_{\Lambda}) & \approx1.81,\quad D_{1}=D_{2}\approx3.51,\quad\left\Vert B_{\OK}\right\Vert _{\infty}\approx4.16,\quad\text{and }\quad\left\Vert B_{\OK}^{-1}\right\Vert _{\infty}=1.
\end{align*}
Consider now again Algorithm \ref{alg:closest-cusp-2} applied to
$\vec{z}=\vec{y}=i\vec{1}\in\Hn$. Since $y_{1}=y_{2}=1$ the first
bounds are given by $a_{1}=a_{2}=D_{0}\approx3.51$ and $P_{\sigma}$
is the polygon bounded by the vertices 
\[
B_{\OK}^{-1}((\pm D_{0},\pm D_{0}))\approx\{(-3.51,0),(0,1.11),(0,-1.11),(3.51,0)\}.
\]
For this particular $\vec{z}$ the preliminary reduction does not
produce any better cusp than $\infty$ and the norm bound is given
by $C\N\left(\vec{y}\right)^{-1/2}=2$ . The domain $P_{\sigma}$
and its pre-image together with the embedded points of $\OK$ and
$\mathbb{Z}^{2}$ are shown in Figure \ref{fig:K2.z1.domains}. We
see that there are only $3$ possibilities for $\sigma$: $-1$, $0$
and $1$ and these result in $10$ candidate cusps. Comparing all
these we see that the cusps $\infty$ and $0$ are both closest with
distance $\Delta(\vec{z},\infty)=\Delta(\vec{z},0)=1$. 

Changing to the point $\vec{z}=i\vec{\frac{1}{2}},$ the preliminary
search finds a tentative closest cusp $0$ with a distance of $1/2$
so we can use the algorithm with $d=1/2,$ which results in the same
bounds for $\sigma$ as before and we find three candidate cusps $0,$
$1$ and $-1$ with the cusp $0$ being the unique closest cusp, with
distance $\Delta\left(\vec{z},0\right)=1/2$. 

To demonstrate the the algorithm works for other cusps than infinity,
consider the point $\vec{z}=(2.58+0.5i,0.5+0.5i)$. The preliminary
search gives only the potential closest cusp $\infty$ so we will
apply the algorithm with $d=2$ and the norm bound $\left|\sigma_{1}\sigma_{2}\right|\le8$.
We find $13$ candidates for $\sigma$ and $35$ distinct candidate
cusps, from which we find that $\mu=(\sqrt{10}:\sqrt{10}+2)$ is the
closest, with distance $\approx1.59$. It is not hard to check that
$\mu$ is equivalent to $\lambda_{2}$ under the element
\[
\left(\begin{array}{cc}
-5 & -2\sqrt{10}+9\\
-2\sqrt{10}+1 & 4\sqrt{10}-10
\end{array}\right)\in\Gamma_{K}.
\]
Applying the complete reduction map to $\vec{z}$ gives $\vec{w}=B\vec{z}$
with $\vec{w}\approx(-0.669+0.036i,0.709+0.004i)$ and 
\[
B=\left(\begin{array}{cc}
-2\sqrt{10}-9 & 9\\
-4\sqrt{10}-9 & 4\sqrt{10}
\end{array}\right)\in\Gamma_{K}.
\]
\end{example}

\begin{example}
To demonstrate that the method works also in degree $3$, consider
$K_{3}=\mathbb{Q}\left(\alpha\right)$ where $\alpha$ has minimal
polynomial $\alpha^{3}-\alpha^{2}-2\alpha+1$. This field has degree
$3$, class number $1$, discriminant $49$, fundamental units $\varepsilon_{1}=2-\alpha^{2}$
and $\varepsilon_{2}=\alpha^{2}-1$ and $\OK$  has an
integral basis 
\[
\beta_{1}=1,\ \beta_{2}=\alpha,\ \beta_{3}=\alpha^{2}-2.
\]
The real embeddings of $\alpha$ are approximately $(-1.247,0.445,1.802)$
and we find the relevant numerical bounds to be: 

\begin{align*}
\vec{r}(B_{\Lambda}) & \approx(1.40,0.810,1.03),\quad\vec{D}\approx(2.01,1.499,1.674),\\
\left\Vert B_{\OK}\right\Vert _{\infty} & \approx4.05,\text{ and }\left\Vert B_{\OK}^{-1}\right\Vert _{\infty}=1.
\end{align*}
Let $\vec{z}=i\vec{1}$ and apply Algorithm \ref{alg:closest-cusp-2}
to find closest cusps. In the first step we find $5$ candidates for
$\sigma$. See Figure \ref{fig:K3.z1.domains}, which shows the polyhedron
together with the surfaces $\N\left(\sigma\right)=1$. The two points
which do not satisfy the norm bound are drawn in lighter gray, the
others in black. In the end we find $8$ candidates for the closest
cusp and we find (as usual) that the cusps $\infty$ and $0$ are
both closest with a distance of $1$. 
\begin{figure}
\begin{minipage}[t]{0.49\linewidth}
\centering\resizebox{\linewidth}{!}{\includegraphics[scale=0.25]{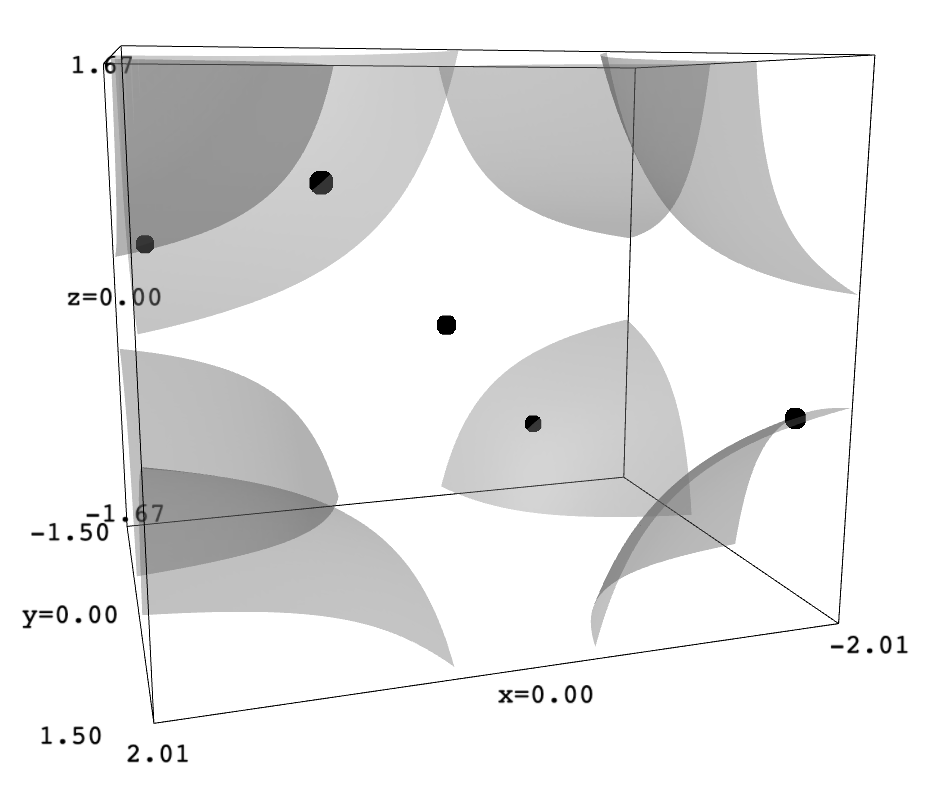}}

\end{minipage}     
\begin{minipage}[t]{0.49\linewidth}

\resizebox{\linewidth}{!}{\includegraphics[scale=0.25]{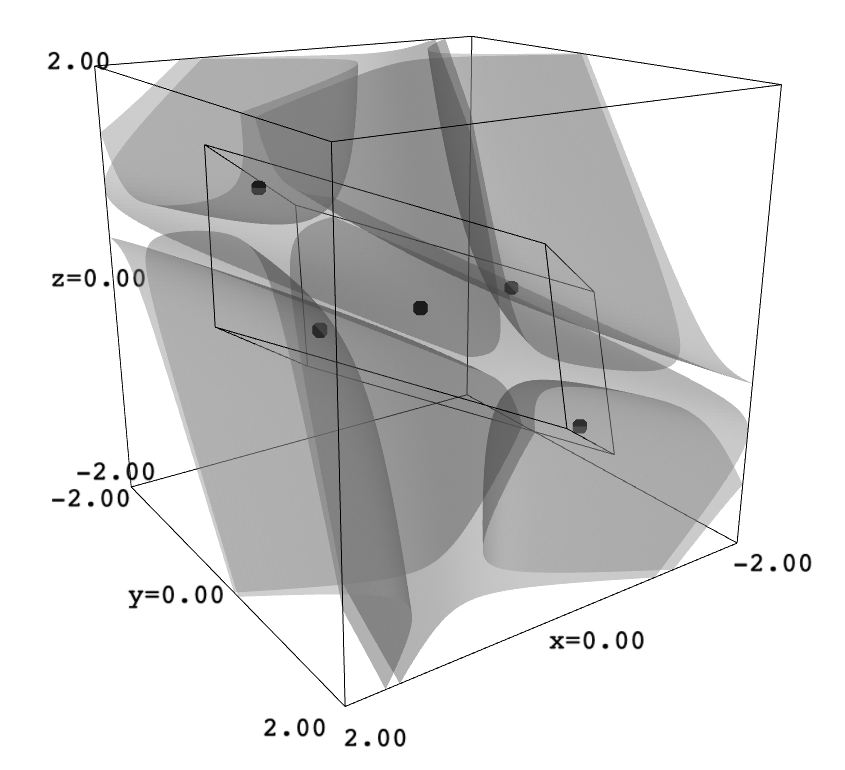}}

%\label{Figure2}         
%\caption{$P_{\sigma}$ for $\mathbf{z}=i\mathbf{1}$}
\end{minipage}

\caption{\label{fig:K3.z1.domains} $B_{\mathcal{O}_{K_3}}(P_{\sigma})$  and $P_{\sigma}$ for $\mathbf{z}=i\mathbf{1}$, $K_3=\mathbb{Q}(\alpha)$ with $\alpha^{3}-\alpha^{2}-2\alpha+1$}
\end{figure}
\end{example}

\begin{example}
Just to give an idea of how it works in a more complicated example,
consider $K=\mathbb{Q}\left(\alpha\right)$ where $\alpha$ has a
minimal polynomial $x^{3}-36x-1$. Then $K$ has discriminant $20733$,
class number $5$ and its label in the LMFDB is 3.3.20733.1. The fundamental
units are $\varepsilon_{1}=-\alpha$ and $\varepsilon_{2}=-\alpha-6$
and  $\OK$ has an integral basis 
\[
\beta_{1}=1,\ \beta_{2}=\alpha,\ \beta_{3}=\frac{1}{3}\left(\alpha^{2}+\alpha-23\right).
\]
The real embeddings of $\alpha$ are approximately $(-5.986,-0.028,6.014)$
and we find the relevant numerical bounds to be: 

\begin{align*}
\left\Vert B_{\Lambda}\right\Vert _{\infty} & \approx6.06,\quad D_{0}\approx52.22,\quad\left\Vert B_{\OK}\right\Vert _{\infty}\approx13.41,\quad\text{and }\quad\left\Vert B_{\OK}^{-1}\right\Vert _{\infty}=1.
\end{align*}
Using Algorithm \ref{alg:closest-cusp-2} with $\vec{z}=i\vec{1}$
we find $9$ candidates for $\sigma$, in total $3396$ candidate
cusps and as usual the cusps $0$ and $\infty$ are both closest with
distance $1$.
\end{example}

\section{Conclusion and Suggestions for Future Work}

It is important to stress again that the raison d'être and main result
of this paper is that, as far as the author is aware, we are for the
first time presenting a reduction algorithm for Hilbert modular groups
which applies to any totally real number field (in theory at least)
and can be proven to return a reduced point and terminates in polynomially
bounded time for a fixed field and varying $\vec{z}$. 

\subsection{Motivation and future applications}

Our interest in reduction theory for Hilbert modular groups stems
from two different problems. The first problem is regarding dimension
formulas for vector-valued Hilbert modular forms. This is part of
ongoing work joint with Skoruppa and Boylan, cf. e.g. \cite{SkoStr17}
and \cite{MR3973298}. One of the necessary ingredients for dimension
formulas is the number of elliptic fixed points, and in the vector-valued
case it is also necessary to know the corresponding stabilizers. The
number of elliptic fixed points is well known for quadratic fields
but for higher degrees this is a hard problem for which a computational
approach currently seems to be the only option. While there are many
computational approaches, both algebraic and analytic, at some point
they generally require some form of reduction to produce representative
elements. 

The second problem is the computation
of non-holomorphic Hilbert modular forms. One of the key ingredients
in the so-called automorphy (or Hejhal's) method for computing Maa\ss cusp forms 
on Hecke triangle groups and subgroups of the modular groups
is the existence of an efficient reduction algorithm. While many parts
of this algorithm need to be modified to work over fields other than
$\mathbb{Q}$, the main obstacle so far has been the lack of a general
reduction algorithm. With the existence of the current algorithm the
hope is that a computational approach to non-holomorphic Hilbert modular
forms is finally within reach. 

From an algorithmic perspective it is clear the most important improvment
would be to find a better bound for the embeddings or coordinates
of $\rho$. While we believe that most of the bounds are close to
optimal in the general setting it might be possible to hard-code the
case of, say, a quadratic field, more efficiently. 

\def\polhk#1{\setbox0=\hbox{#1}{\ooalign{\hidewidth
  \lower1.5ex\hbox{`}\hidewidth\crcr\unhbox0}}}

\end{document}